\documentclass[12pt]{article}
\usepackage{amsmath,amsthm,amsfonts,amssymb,latexsym}

\author{Constantin-Nicolae Beli} 
\title{Universal integral quadratic forms over dyadic local fields} 
\date{}

\newtheorem{theorem}{Theorem}[section]
\newtheorem{proposition}[theorem]{Proposition}
\newtheorem{lemma}[theorem]{Lemma}
\newtheorem{definition}{Definition}
\newtheorem{corollary}[theorem]{Corollary}

\newtheorem{bof}[theorem]{}
\newtheorem{teorema}{Theorem}

\def\qed{\mbox{$\Box$}\vspace{\baselineskip}}
\def\pf{$Proof.$ } 
\def\bco{\begin{corollary}} \def\eco{\end{corollary}} 
\def\bdf{\begin{definition}} \def\edf{\end{definition}} 
\def\btm{\begin{theorem}} \def\etm{\end{theorem}} 
\def\bpr{\begin{proposition}} \def\epr{\end{proposition}}  
\def\blm{\begin{lemma}} \def\elm{\end{lemma}} 
\def\bff{\begin{bof}\rm} \def\eff{\end{bof}}
\def\btr{\begin{teorema}} \def\etr{\end{teorema}}

 \def\ZZ{{\mathbb Z}} \def\QQ{{\mathbb Q}}
 
\def\fff{\dot F} \def\ffs{\dot F^2} \def\e{\varepsilon}
\def\oo{{\mathcal O}} \def\ooo{\oo^\times}
\def\oos{\oo^{\times 2}} \def\www{{\mathfrak w}} \def\FFF{{\mathfrak f}} 
\def\p{{\mathfrak p}} \def\sss{{\mathfrak s}} \def\nnn{{\mathfrak n}}
\def\rep{{\to\!\!\! -}}
\DeclareMathOperator\ord{ord} \DeclareMathOperator\rank{rank}

\newtheorem{lem}{Lema}
{\begin{lem}\rm}{\end{lem}}
\begin{document}
\textwidth 6.5 truein 
\oddsidemargin 0 truein 
\evensidemargin -0.50 truein
\topmargin -.5 truein
\textheight 9 in

\maketitle

\begin{abstract}
We give necessary and sufficient conditions for an integral quadratic
form over a dyadic local field to be universal.
\end{abstract}





\section{Introduction}

Let $F$ be a local non-archimedian field of characteristic zero. Let
$\oo$ be the ring of integers and let $\p$ be the prime ideal of
$F$. The group of units of $\oo$ is $\ooo =\oo\setminus\p$. We have
$\p =\pi\oo$, where $\pi$ is a prime element of $\oo$. If $\mathfrak
a$ is an ideal of $F$ then we define its order $\ord\mathfrak
a\in\ZZ\cup\{\infty\}$ by $\ord\mathfrak a=R$ if $\mathfrak a=\p^R$
for some $R\in\ZZ$ and $\ord\mathfrak a=\infty$ if $\mathfrak a=0$. If
$a\in F$ we denote by $\ord a=\ord a\oo$, i.e. $\ord a$ is the value
of $a$. We have $\ord a=R$ if $a=\pi^R\varepsilon$ with
$\varepsilon\in\ooo$ and $\ord a=\infty$ if $a=0$.

We denote by $(\cdot,\cdot )_\p\to\fff/\fff^2\times\fff/\fff^2\to\{\pm
1\}$ the Hilbert symbol.

All quadratic spaces and lattices in this paper will be assumed to be
non-degenerate.

If $V$ is a quadratic space and $x_1,\ldots,x_n$ is an orthogonal
basis with $Q(x_i)=a_i$, then we say that $V\cong [a_1,\ldots,a_n]$
relative to the orthogonal basis $x_1,\ldots,x_n$. For the quadratic
lattice $L=\oo x_1\perp\cdots\perp\oo x_n$ we write $L\cong\langle
a_1,\ldots,a_n\rangle$. 

Recall that if $b\in\fff$ then $b$ is represented by $[a_1,a_2]$ iff
$(a_1b,-a_1a_2)_p=1$ and it is represented by $[a_1,a_2,a_3]$ iff
$b\notin -a_1a_2a_3\fff^2$ or $[a_1,a_2,a_3]$ is isotropic. We also
have that $[a_1,a_2,a_3]$ is isotropic iff $-a_1$ is represented by
$[a_2,a_3]$, which is equivalent to $(-a_1a_2,-a_2a_3)_p =1$.

If $V,W$ are two quadratic spaces, we denote by $W\rep V$ the
fact that $V$ represents $W$. Similarly for lattices.

If $L$ is a quadratic lattice, with $FL=V$, and $Q:V\to F$ is the
corresponding quadratic form, then we say that $L$ is integral if
$Q(L)\subseteq\oo$ and we say that it is universal if $Q(L)=\oo$. In
[XZ] the authors gave necessary and sufficient conditions for a
quadratic lattice to be universal in the case when $F$ is
non-dyadic. In the more complicated dyadic case they solved the same
problem, but only for binary and ternary lattices.

In this paper we completely solve this problem for dyadic quadratic
lattices in arbitrary dimensions. Unlike in [XZ], where the quadratic
lattices are described in terms of Jordan compositions, here we use
BONGs (bases of norm generators), which we introduced in [B1]. Since
the BONGs are not widely known and used, we now give a brief
review. A summary of the results from [B1] we use here can be found in
[B3, \S1].

From now on $F$ is a dyadic field, i.e. a finite extensions of
$\QQ_2$. We denote by $e$ the ramification index of the extension
$F/\QQ_2$, i.e. $e=\ord 2$.

\subsection{The map
$d:\fff/\fff^2\to\{0,1,3,5,\ldots,2e-1,2e,\infty\}$}

The quadratic defect, introduced in [OM, \S63A], of an element
$\varepsilon\in F$ is the ideal $\mathfrak d(a)=\bigcap_{x\in
F}(a-x^2)\oo$. We denote by $\Delta =1-4\rho$ a fixed element with
$\mathfrak d(\Delta )=4\oo$.

In [B1, \S1] we introduced the order of the relative quadratic defect
$d:\fff/\fff^2\to\ZZ\cup\{\infty\}$, $d(a)=\ord a^{-1}\mathfrak
d(a)$. Let $a=\pi^R\varepsilon$, with $\varepsilon\in\ooo$. If $R$ is
even $d(a)=d(\varepsilon )=\ord\mathfrak d(\varepsilon )\in\{
1,3,5,\ldots,2e-1,2e,\infty\}$. If $R$ is odd then $d(a)=0$.

We have $d(a)=0$ iff $\ord a$ is odd, $d(a)\geq 1$ iff $\ord a$ is
even, $d(a)=2e$ iff $a\in\Delta\fff^2$ and $d(a)=\infty$ iff
$a\in\fff^2$.

The map $d$ has the folowing properties:

\noindent (1) $d(ab)\geq\min\{ d(a),d(b)\}$ $\forall a,b\in\fff$.

\noindent (2) If $d(a)+d(b)>2e$ then $(a,b)_\p =1$.

\noindent (3) If $a\in\fff\setminus\fff^2$ then there is $b\in\fff$ with
$d(b)=2e-d(a)$ such that $(a,b)_\p =-1$. Moreover, if $d(a)<2e$ then
we can choose $b\in\ooo$. 

(For the last statement note that if $d(a)<2e$ then $d(b)=2e-d(a)>0$
so $b\in\ooo\fff^2$. Since both $d(b)$ and $(a,b)_\p$ depend only on
$b$ modulo $\fff^2$, we may assume that $b\in\ooo$.)


\subsection{BONGs and good BONGs}

Let $V$ be a quadratic space over $F$, with the quadratic form $Q:V\to
F$ and the corresponding bilinear symmetric form $B:V\times V\to F$,
$B(x,y)=\frac 12(Q(x+y)-Q(x)-Q(y))$.

Let now $L$ be a lattice over $V$. The norm $\nnn L$ of $L$ is the
fractionary ideal generated by $Q(L)$ and the scale $\sss L$ of $L$ is
the fractionary ideal $B(L,L)$. 

{\bf Bases of norm generators (BONGs)}

The bases of norm generators (BONGs), introduced in [B1, \S2], were
defined recursively as follows. A norm generator of $L$ is and element
$x\in L$ such that $Q(x)\oo =\nnn L$. A basis of norm generator (BONG)
of $L$ is an orthogonal basis $x_1,\ldots,x_n$ of $V=FL$ such that
$x_1$ is a norm generator for $L$ and $x_2,\ldots,x_n$ are a BONG for
$pr_{x_1^\perp}L$. (Here $pr_{x_1^\perp}:V\to x_1^\perp$ is the
projection on the othogonal complement $x_1^\perp$ of $x_1$.)

A BONG uniquely determines a lattice. We write $L=\prec
x_1,\ldots,x_n\succ$ to denote the fact that $x_1,\ldots,x_n$ is a
BONG for $L$ and we say that $L\cong\prec a_1,\ldots,a_n\succ$
relative to the BONG $x_1,\ldots,x_n$ if $Q(x_i)=a_i$.

{\bf The binary case}

If $n=2$ then an orthogonal basis $x_1,x_2$ of $V$ with $Q(x_i)=a_i$
is the BONG of a lattice iff $a_2/a_1\in\mathcal A$,
where $\mathcal A\subset\fff/\oos$, $\mathcal A=\{ a\in\fff\,\mid\,
a\in\frac 14\oo,\,\mathfrak d(-a)\subseteq\oo\}$.

If $a\in\fff$ with $\ord a=R$, then $a\in\mathcal A$ iff $R+d(-a)\geq
0$ and $R\geq -2e$. Hence if $\ord a_i=R_i$ then $a_2/a_1\in\mathcal
A$ iff $R_2-R_1\geq -2e$ and $R_2-R_1+d(-a_1a_2)\geq 0$. (See [B1,
Lemmas 3.5 and 3.6].)

If $R_2>R_1$ then we have the Jordan splitting $L=\oo x_1\perp\oo x_2$
and the scales of the Jordan components are $\nnn\oo x_1=\p^{R_1}$ and
$\nnn\oo x_2=\p^{R_2}$.

If $R_2\leq R_1$ then $L$ is $\p^{(R_1+R_2)/2}$-modular with $\nnn
L=\p^{R_1}$. In particular, if $R_2-R_1=-2e$, then $\sss
L=\p^{(R_1+R_2)/2}=\p^{R_1-e}=\frac 12\p^{R_1}=\frac 12\nnn L$. Hence
$L\cong\frac 12\p^{R_1}A(0,0)$ or $\frac 12\p^{R_1}A(2,2\rho )$.

If $a\in{\mathcal A}$, with $\ord a=R$, then $g(a)\leq\ooo/\oos$ is
defined by $g(a)=\ooo$ if $R=-2e$, $g(a)=\oos$ if $R>2e$ and
$$g(a)=\begin{cases}(1+\p^{R/2+e})\oos & d(-a)>e-R/2\\
(1+\p^{R+d(-a)})\oos\cap{\rm N}(-a) & d(-a)\leq e-R/2\end{cases}$$
if $-2e<R\leq 2e$.

Then if $L\cong\prec a_1,a_2\succ$ and $\eta\in\ooo$, we have $L\cong
L^\eta$, i.e. $\prec a_1,a_2\succ\cong\prec\eta a_1,\eta a_2\succ$,
iff $\eta\in g(a_2/a_1)$. (See [B1, Lemma 3.11].)

{\bf Good BONGs}

A BONG $x_1,\ldots,x_n$ of $L$ is called good if $\ord Q(x_i)\leq\ord
Q(x_{i+2})$ for $1\leq i\leq n-2$.

If $x_1,\ldots,x_n$ is an orthogonal basis of $V$, with $Q(x_i)=a_i$
and $\ord a_i=R_i$ then $x_1,\ldots,x_n$ is the good BONG of a lattice
$L$ with $FL=V$ iff $R_i\leq R_{i+2}$ for $1\leq i\leq n-2$ and
$a_{i+1}/a_i\in\mathcal A$ for $1\leq i\leq n-1$. The second condition
writes as $R_{i+1}-R_i\geq -2e$ and $R_{i+1}-R_i+d(-a_ia_{i+1})\geq
0$. (See [B1, Lemma 4.3(ii)].)

In particular, if $R_{i+1}-R_i$ is odd then $\ord
a_ia_{i+1}=R_i+R_{i+1}$ is odd so
$R_{i+1}-R_i=R_{i+1}-R_i+d(-a_ia_{i+1})\geq 0$. Thus $R_{i+1}-R_i$
cannot be odd and negative.

If $R_{i+1}-R_i=-2e$ then $R_{i+1}-R_i+d(-a_ia_{i+1})\geq 0$ implies
$d(-a_ia_{i+1})\geq 2e$ so $-a_1a_2\in\fff^2$ or $\Delta\fff^2$,
corresponding to $d(-a_1a_2)=\infty$ or $2e$, accordingly.

Every quadratic lattice has a good BONG. Good BONGs can be obtained
with the help of the so-called maxinmal norm splittings. (See [B1,
Lemmas 4.3(iii) and 4.6] and [B3, \S7].)

{\bf Similarities with orthogonal bases}

Unlike in the non-dyadic case, in the dyadic case lattices usually
don't have orthogonal bases. The BONGS, especially the good BONGs, are
a good substitute, as the they preserve many of the properties of the
orthogonal bases.

Suppose that $x_1,\ldots,x_n$ are an orthogonal basis of a quadratic
space, with $Q(x_i)=a_i$ and $\ord a_i=R_i$.

If $x_1,\ldots,x_n$ is a good BONG for $L$ then $L=\prec x_1,\ldots
x_k\succ\perp\prec x_{k+1},\ldots,x_n\succ$ holds iff $R_k\leq
R_{k+1}$. Equivalently, $\prec a_1,\ldots,a_n\succ\cong\prec a_1,\ldots
a_k\succ\perp\prec a_{k+1},\cdots,a_n\succ$ iff $R_k\leq R_{k+1}$.

Conversely, if $x_1,\ldots,x_k$ and $x_{k+1},\ldots,x_n$ are good
BONGs for the lattices $L'$ and $L''$ then $x_1,\ldots,x_n$ is a good
BONG for $L'\perp L''$ iff $R_k\leq R_{k+1}$, $R_{k-1}\leq R_{k+1}$
and $R_k\leq R_{k+2}$. (If $k=1$ we ignore $R_{k-1}\leq R_{k+1}$; if
$k=n-1$ we ignore $R_k\leq R_{k+2}$.)

If $x_1,\ldots,x_n$ is a good BONG and for some $1\leq k\leq l\leq n$
$y_k,\ldots,y_l$ is another good BONG for $\prec y_k,\ldots,y_l\succ$,
then $x_1,\ldots,x_{k-1},y_k,\ldots,y_l,x_{l+1}\ldots,x_n$ is a good
BONG for $L$. Consequently, if $L\cong\prec a_1,\ldots,a_n\succ$ and
$\prec a_k,\ldots,a_l\succ\cong\prec b_k,\ldots,b_l\succ$ relative to
good BONGs, then ${L\cong\prec
a_1,\ldots,a_{k-1},b_k,\ldots,b_l,b_{l+1}\ldots,b_n\succ}$ relative to
a good BONG.

In particular, if $\eta\in g(a_{k+1}/a_k)$ then $\prec
a_k,a_{k+1}\succ\cong\prec\eta a_k,\eta a_{k+1}\succ$ so\\ ${\prec
a_1,\ldots,a_n\succ\cong\prec a_1,\ldots,a_{k-1},\eta a_k,\eta
a_{k+1},a_{k+2}\ldots a_n\succ}$.

If $L\cong\prec a_1,\ldots,a_n\succ$ relative to the good BONG
$x_1,\ldots,x_n$ then $L^\sharp\cong\prec
a_n^{-1},\ldots,a_1^{-1}\succ$ relative to the good BONG
$x_n^\sharp,\ldots,x_1^\sharp$, where $x^\sharp :=Q(x)^{-1}x$ for
every $x\in V$ with $Q(x)\neq 0$. 

(See [B1, Lemmas 4.8 and 4.9, Corollary 4.4].)

\subsection{The invariants $R_i(L)$ and $\alpha_i(L)$ and the
clasification theorem}

Suppose now that $L\cong\prec a_1,\ldots,a_n\succ$ relative to some
good BONG and let $R_i=\ord a_i$. In [B3, \S2] we defined, for every
$1\leq i\leq n-1$, the number $\alpha_i$ as the minimum of the set

\begin{multline*}
\{ (R_{i+1}-R_i)/2+e\}\cup\{ R_{i+1}-R_j+d(-a_ja_{j+1})\,\mid\,
1\leq j\leq i\}\\
\cup\{ R_{j+1}-R_i+d(-a_ja_{j+1})\,\mid\,
i\leq j\leq n-1\}.
\end{multline*}

The numbers $R_i$, with $1\leq i\leq n$, and $\alpha_i$, with $1\leq
i\leq n-1$, are invariants of the lattice $L$, so we denote
them by $R_i(L)$ and $\alpha_i(L)$. If $L$ has a Jordan
decomposition $L=L_1\perp\cdots\perp L_t$, then the numbers
$R_i=R_i(L)$ are in one-to-one correpondence with $t$, $\rank L_i$,
$\sss L_i$ and $\nnn L^{\sss L_i}$ with $1\leq i\leq t$. In
particular, $\nnn L=\p^{R_1}$ and $\sss L=\p^{\min\{
R_1,(R_1+R_2)/2\}}$. The numbers $\alpha_i=\alpha_i(L)$ are in
one-to-one correspondence with the invariants $\www_i=\www L^{\sss
L_i}$, with $1\leq i\leq t$, and $\mathfrak f_i$, with $1\leq i\leq
t-1$, of $L$. (See [B1, Lemma 4.7] and [B3, Lemmas 2.13(i), 2.15 and
2.16].) 
\medskip

Here are some properties of the invariants $\alpha_i$, which appear in
[B3, Lemmas 2.2 and 2.7, Corollaries 2.8 and 2.9, Remark 2.6].

(1) The sequence $R_i+\alpha_i$ is increasing and the sequence
$-R_{i+1}+\alpha_i$ is decreasing.

(2) $\alpha_i\geq 0$, with equality iff $R_{i+1}-R_i=-2e$.

(3) If $R_{i+1}-R_i\leq 2e$ then $\alpha_i\geq R_{i+1}-R_i$ with
equality iff $R_{i+1}-R_i=2e$ or it is odd.

(4) If $R_{i+1}-R_i\in\{ -2e,2-2e,2e-2\}$ or $R_{i+1}-R_i\geq 2e$ then
$\alpha_i=(R_{i+1}-R_i)/2+e$.

(5) $\alpha_i$ is $<2e$ $=2e$ or $>2e$ iff $R_{i+1}-R_i$ is so.

(6) $\alpha_i$ is an odd integer unless $\alpha_i=(R_{i+1}-R_i)/2+e$.

(7) $\alpha_i$ is an integer unless $R_{i+1}-R_i$ is odd and $>2e$.

(8) $\alpha_i\in ([0,2e]\cap\ZZ )\cup ((2e,\infty )\cap\frac 12\ZZ
))$. 

(9) $\alpha_i=\min\{ (R_{i+1}-R_i)/2+e,\,
R_{i+1}-R_i+d(-a_ia_{i+1}),\, R_{i+1}-R_i+\alpha_{i-1},\,
R_{i+1}-R_i+\alpha_{i+1}\}$. (If $i=1$ we ignore
$R_{i+1}-R_i+\alpha_{i-1}$, as $\alpha_0$ is not defined. If $i=n-1$
then we ignore $R_{i+1}-R_i+\alpha_{i+1}$, as $\alpha_n$ is not
defined.)

(10) $\alpha_i$ are invariant to scalling.

(11) $R_i(L^\sharp )=-R_{n+1-i}$ and $\alpha_i(L^\sharp )=\alpha_{n-i}(L)$.
\medskip

We now state O'Meara's classification theorem [OM, Theorem 93:28] in
terms of BONGs. This result is [B3, Theorem 3.1].

\btm Let $L,K$ be two quadratic lattice with $FL\cong FK$ and let
$L\cong\prec a_1,\ldots,a_n\succ$ and $K\cong\prec
b_1,\ldots,b_n\succ$ relatice to good BONGs. Let $R_i=R_i(L)$,
$S_i=R_i(K)$, $\alpha_i=\alpha_i(L)$ and $\beta_i=\alpha_i(K)$. Then
$L\cong K$ iff:

(i) $R_i=S_i$ for $1\leq i\leq n$

(ii) $\alpha_i=\beta_i$ for $1\leq i\leq n-1$.

(iii) $d(a_1\cdots a_i\, b_1\cdots b_i)\geq\alpha_i$ for $1\leq i\leq
n-1$.

(iv) $[b_1,\ldots,b_{i-1}]{\rep}[a_1,\ldots,a_i]$ for every
$1<i<n$ such that $\alpha_{i-1}+\alpha_i>2e$.
\etm

\subsection{The invariants $d[\varepsilon a_{i,j}]$ and $d[\varepsilon
a_{1,i}b_{1,j}]$}

For convenience, if $a_1,a_2,\ldots\in\fff$ and $1\leq i\leq j+1$ then
we denote by $a_{i,j}=a_i\cdots a_j$. By convention, $a_{i,i-1}=1$.

If $L\cong\prec a_1,\ldots,a_n\succ$ relative to a good BONG and
$\alpha_i=\alpha_i(L)$, then for every $0\leq i-1\leq j\leq n$ and
$\varepsilon\in\fff$, then we define
$$d[a_{i,j}]=\min\{ d(a_{i,j}),\alpha_{i-1},\alpha_j\}.$$
If $i-1\in\{ 0,n\}$ $\alpha_{i-1}$ is not defined so it is
ignored. Similarly $\alpha_j$ is ignored if $j\in\{ 0,n\}$.)

In particular, since $d[-a_{i,i+1}]=\min\{
d(-a_{i,i+1}),\alpha_{i-1},\alpha_{i+1}\}$, the property (9) of \S1.3
can be written as
$$\alpha_i=\min\{ (R_{i+1}-R_i)/2+e,R_{i+1}-R_i+d[-a_{i,i+1}]\}.$$

If $M\cong\prec a_1,\ldots,a_m\succ$ and $N\cong\prec
b_1,\ldots,b_n\succ$ relative to good BONGs, $\alpha_i=\alpha_i(M)$
and $\beta_i=\alpha_i(N)$, then for every $0\leq i\leq m$, $0\leq
j\leq m$, we define
$$d[\varepsilon a_{1,i}b_{1,j}]=min\{ d(\varepsilon
a_{1,i}b_{1,j}),\alpha_i,\beta_j\}.$$
(If $i\in\{ 0,m\}$ then we ignore $\alpha_i$. If $j\in\{ 0,n\}$ then
we ignore $\beta_i$.)

As a consequence of condition (iii) of Theorem ..., $d[\varepsilon
a_{i,j}]$ and $d[\varepsilon a_{1,i}b_{1,j}]$ are independent of the
choice of the good BONGs. Also $d[\varepsilon a_{i,j}]$ are a
particular case of the expression $d[\varepsilon
  a_{1,i}b_{1,j}]$. Indeed, if we take $M=N=L$, so that $b_i=a_i$ and
$\beta_i=\alpha_i$ then in $\fff/\fff^2$ we have
$a_{1,j}b_{1,i-1}=a_{1,j}a_{1,i-1}=a_{i,j}$ so $d(\varepsilon
a_{1,j}b_{1,i-1})=d(\varepsilon a_{i,j})$. Therefore
$$d[\varepsilon a_{1,j}b_{1,i-1}]=\min\{ d(\varepsilon
a_{1,j}b_{1,i-1}),\alpha_j,\beta_{i-1}\}=\min\{d(\varepsilon
a_{i,j}),\alpha_j,\alpha_{i-1}\}=d(\varepsilon a_{i,j}).$$

The invariants $d[\cdot ]$ satisfy a similar domination principle as
$d(\cdot )$. Namely, if we have a third lattice $K\cong\prec
c_1,\ldots,c_k\succ$ and $\varepsilon,\varepsilon'\in\fff$ then in
$\fff/\fff^2$ we have $(\varepsilon
a_{1,i}b_{1,j})(\varepsilon'b_{1,j}c_{1,k})
=\varepsilon\varepsilon'a_{1,i}c_{1,k}$ and so
$d(\varepsilon\varepsilon'a_{1,i}c_{1,k})\geq\min\{ d(\varepsilon
a_{1,i}b_{1,j}),d(\varepsilon'b_{1,j}c_{1,k})\}$. Similarly, we have
$$d[\varepsilon\varepsilon'a_{1,i}c_{1,k}]\geq\min\{ d[\varepsilon
a_{1,i}b_{1,j}],d[\varepsilon'b_{1,j}c_{1,k}]\}.$$

Note that both $d$ and $\alpha_i$ take nonnegative values so
$d[\varepsilon a_{1,i}b_{1,j}]$ is always nonnegative. If
$\ord\varepsilon a_{1,i}b_{1,j}$ is odd then $d(\ord\varepsilon
a_{1,i}b_{1,j})=0$ and so $d[\ord\varepsilon a_{1,i}b_{1,j}]=0$.

\subsection{The representation theorem}

We now state the representation theorem, which was announced in [B2,
Theorem 4.5].

Let $M,N$ be quadratic lattices, with
$M\cong\prec a_1,\ldots,a_m\succ$ and
$N\cong\prec b_1,\ldots,b_n\succ$ relative to good BONGs and
$m\geq n$. Let $R_i=R_i(M)$, $S_i=R_i(N)$, 
$\alpha_i=\alpha_i(M)$ and $\beta_i=\alpha_i(N)$. If $1\leq
i\leq\min\{ m-1,n\}$, then we define $A_i=A_i(M,N)$ as
$$A_i=\min\{ (R_{i+1}-S_i)/2+e,\, R_{i+1}-S_i+d[-a_{1,i+1}b_{1,i-1}],
R_{i+1}+R_{i+2}-S_{i-1}-S_i+d[a_{1,i+2}b_{1,i-2}]\}.$$
(If $i=1$ or $m-1$, then the term
$R_{i+1}+R_{i+2}-S_{i-1}-S_i+d[a_{1,i+2}b_{1,i-2}]$ is not defined so
it is ignored.)

If $n\leq m-2$ the we assume that $S_{n+1}\gg 0$. Then, formally, we
have
$$S_{n+1}+A_{n+1}=\min\{ (R_{n+2}+S_{n+1})/2+e,\,
R_{n+2}+d[-a_{1,n+2}b_{1,n}],\,
R_{n+2}+R_{n+3}-S_n+d[a_{1,n+3}b_{1,n-1}]\}.$$
Since $(R_{n+2}+S_{n+1})/2+e\to\infty$ as $S_{n+1}\to\infty$, we can
ignore it in the above formula and we define
$$S_{n+1}+A_{n+1}=\min\{ R_{n+2}+d[-a_{1,n+2}b_{1,n}],\,
R_{n+2}+R_{n+3}-S_n+d[a_{1,n+3}b_{1,n-1}]\}.$$
(If $n=m-2$ then $R_{n+2}+R_{n+3}-S_n+d[a_{1,n+3}b_{1,n-1}]$ is not
defined, so we ignore it.)

\btm Assume that $FM\rep FN$. Then $M\rep N$ iff:

(i) For $1\leq i\leq n$ we have $R_i\leq S_i$ or $1<i<m$ and
$R_i+R_{i+1}\leq S_{i-1}+S_i$.

(ii) For $1\leq i\leq\min\{ m-1,n\}$ we have $d[a_{1,i}b_{1,i}]\geq
A_i$.

(iii) For any $1<i\leq\min\{ m-1,n+1\}$ such that $R_{i+1}>S_{i-1}$
and $A_{i-1}+A_i>2e+R_i-S_i$ we have $[b_1,\ldots,b_{i-1}]\rep
[a_1,\ldots,a_i]$.

(iv) For any $1<i\leq\min\{ m-2,n+1\}$ such that $R_{i+1}\leq
S_{i-1}$, $R_{i+2}\leq S_i$ and $R_{i+2}-S_{i-1}>2e$ we have
$[b_1,\ldots,b_{i-1}]\rep [a_1,\ldots,a_{i+1}]$. (If $i=n+1$ we
ignore the condition $R_{i+2}\leq S_i$.)
\medskip

Note that if $n\leq m-2$ and $i=n+1$ then $S_{n+1}$ and $A_{n+1}$ are
not defined, but $S_{n+1}+A_{n+1}$ is. Thus the condition
$A_n+A_{n+1}>2e+R_{n+1}-S_{n+1}$ from (iii) should be read as
$A_n+(S_{n+1}+A_{n+1})>2e+R_{n+1}$.
\etm

{\bf Remarks}

{\bf 1.} We have $R_i-S_{i-1}+d[-a_{1,i}b_{1,i-2}]\geq A_{i-1}$ and
$R_{i+1}-S_i+d[-a_{1,i+1}b_{1,i-1}]\geq A_i$, by the definition of
$A_i$. So if $A_{i-1}+A_i>2e+R_i-S_i$ then
$R_i-S_{i-1}+d[-a_{1,i}b_{1,i-2}]+R_{i+1}-S_i+d[-a_{1,i+1}b_{1,i-1}]
>2e+R_i-S_i$
so
$d[-a_{1,i}b_{1,i-2}]+d[-a_{1,i+1}b_{1,i-1}]>2e+S_{i-1}-R_{i+1}$. But,
by [B4, Lemma 3.16], if $M$ and $N$ satisfy the conditions (i) and
(ii) of Theorem 1.2 and $R_{i+1}>S_{i-1}$, then
$A_{i-1}+A_i>2e+R_i-S_i$ is equivalent to
$d[-a_{1,i}b_{1,i-2}]+d[-a_{1,i+1}b_{1,i-1}]>2e+S_{i-1}-R_{i+1}$. Hence
the condition (iii) of Theorem 1.2 can be replaced by:
\medskip

{\em (iii') For any $1<i\leq\min\{ m-1,n+1\}$ such that $R_{i+1}>S_{i-1}$
and $d[-a_{1,i}b_{1,i-2}]+d[-a_{1,i+1}b_{1,i-1}]>2e+S_{i-1}-R_{i+1}$
we have $[b_1,\ldots,b_{i-1}]\rep [a_1,\ldots,a_i]$.}
\bigskip

{\bf 2.} In some cases conditions (ii) and (iii) (or its equivalent
(iii')) of Theorem 1.2 need not be verified. Namely, an index $1\leq
i\leq\min\{ m,n+1\}$ is called {\em essential} if $R_{i+1}>S_{i-1}$
and $R_{i+1}+R_{i+2}>S_{i-2}+S_{i-1}$. (The inequalities that do not make
sense because $R_{i+1}$, $R_{i+2}$, $S_{i-2}$ or $S_{i-1}$ is not
defined are ignored. Then condition (ii) is vacuous at an index $i$ if
both $i$ and $i+1$ are not essentian and condition (iii) is vacuous at
an index $i$ if $i$ is not essential.
\bigskip

{\bf 3.} Condition (iv) of Theorem 1.2 can be replaced by a stronger
condition, where the inequalities $R_{i+1}\leq S_{i-1}$ and
$R_{i+2}\leq S_i$ are ignored.

We have an even more general reasult. If $N\rep M$ and $R_l-S_j>2e$
for some $1\leq l\leq m$, $1\leq j\leq n$, then $[b_1,\ldots,b_j]\rep
[a_1,\ldots,a_{l-1}]$. In all cases but the one described in condition
(iv), this follows from conditions (i), (ii) and (iii).

In particular, if $R_{i+1}-S_i>2e$ then $[b_1,\ldots,b_i]\cong
[a_1,\ldots,a_i]$. By taking determinants we get that in $\fff/\fff^2$
we have $a_{1,i}=b_{1,i}$ or, equivalently
$d(a_{1,i}b_{1,i})=\infty$.

In fact, we have an even stronger result. If $N\rep M$ and
$R_l-S_j>2e$ then\\ $\prec b_1,\ldots,b_j\succ\rep\prec
a_1,\ldots,a_{l-1}\succ$.

\section{The main result}

\btm Let $M$ be an integral quadratic lattice with $M\cong\prec
a_1,\ldots,a_m\succ$ relative to a good BONG and let $R_i=R_i(M)$ for
$1\leq i\leq m$, $\alpha_i=\alpha_i(M)$ for $1\leq i\leq m-1$.

Then $M$ is universal if and only if $m\geq 2$, $R_1=0$ and we have
one of the cases below:
\medskip

I (a) $\alpha_1=0$ or, equivalently, $R_2=-2e$.

~~(b) If $m=2$ or $R_3>1$, then $[a_1,a_2]$ is isotropic.

~~(c) If $m\geq 3$, $R_3=1$ and either $m=3$ or $R_4>2e+1$,
then $[a_1,a_2]$ is isotropic.

\medskip

II (a) $m\geq 3$ and $\alpha_1=1$.

~~(b) If $R_2=1$ or $R_3>1$, then $m\geq 4$ and $\alpha_3\leq
2(e-[\frac{R_3-R_2}2])-1$.

~~(c) If $R_2\leq 0$, $R_3\leq 1$ and either $m=3$ or
$R_4-R_3>2e$, then $[a_1,a_2,a_3]$ is isotropic.
\etm

\medskip

{\bf Remark} By \S1.3, property (2), we have $\alpha_1=0$
iff $R_2-R_1=-2e$. So if $R_1=0$ then $\alpha_1=0$ is equivalent to
$R_2=-2e$.
\bigskip

The condition that $M$ is integral, i.e. that $Q(M)\subseteq\oo$, is
equivalent to $\nnn M\subseteq\oo$. But $\nnn L=\p^{R_1}$ so we have:

\blm $M$ is integral iff $R_1\geq 0$.
\elm

As noted in [XZ, Lemma 2.2], an integral lattice $M$ is universal iff
it represents all elements of $\ooo\cup\pi\ooo$, i.e. the elements
$b\in\fff$ with $\ord b=0$ or $1$. In terms of Theorem 1.2, this means
that $M$ represents every unary lattice $N\cong\prec b_1\succ$ with
$S_1=\ord b_1\in\{ 0,1\}$. 

Then Theorem 1.2 in the case $n=1$ implies that:

\blm The integral lattice $M$ is universal iff for every $N=\prec
b_1\succ$ with $S_1\in\{ 0,1\}$ we have $FN\rep FM$ and:

(i) $R_1\leq S_1$.

(ii) $d[a_1b_1]\geq A_1$.

(iii') If $m\geq 3$, $R_3>S_1$ and
$d[-a_{1,2}]+d[-a_{1,3}b_1]>2e+S_1-R_3$ then $[b_1]\rep  [a_1,a_2]$.

(iv) If $m\geq 4$, $R_3\leq S_1$ and $R_4-S_1>2e$ then $[b_1]\rep
[a_1,a_2,a_3]$.
\elm

Note that for (iii') we used the fact that
$d[-a_{1,2}b_{1,0}]=d[-a_{1,2}]$. As a consequence of this, we also
have $A_1=\min\{ (R_2-S_1)/2+e,\, R_2-S_1+d[-a_{1,2}]\}$.

\blm The condition that $FN\rep FM$ holds for all $N$ iff $FM$
is universal.
\elm
\pf The condition that $FM$ represents every $FN\cong [b_1]$ with
$\ord b_1\in\{ 0,1\}$ is equivalent to $b_1\rep FM$ for every
$b_1\in\fff$, i.e. to $FM$ being universal. \qed

\blm If $M$ is integral, the condition (i) of Lemma 2.3 holds for all
$N$ iff $R_1=0$.
\elm
\pf Since $M$ is integral, we have $R_1\geq 0$. Condition (i) holds
for every $N$ iff $R_1\leq S_1$ for $S_1\in\{ 0,1\}$. Hence the
conclusion. \qed

\blm If $R_1\leq S_1$ then $\alpha_1\geq A_1$, with equality when
$R_1=S_1$.

Consequently, if $R_1\leq S_1$ then $d[a_1b_1]\geq A_1$ is equivalent
to $d(a_1b_1)\geq A_1$.
\elm
\pf We have $A_1=\min\{ (R_2-S_1)/2+e,\, R_2-S_1+d[-a_{1,2}]\}$ and,
by \S1.4, we also have  $\alpha_1=\min\{ (R_2-R_1)/2+e,\,
R_2-R_1+d[-a_{1,2}]\}$. So if $R_1\leq S_1$ then $\alpha_1\geq A_1$
and if $R_1=S_1$ then $\alpha_1=A_1$.

Assume now that $R_1\leq S_1$ so $\alpha_1\geq A_1$. Then, since
$d[a_1b_1]=\min\{ d(a_1b_1),\alpha_1\}$, we have $d[a_1b_1]\geq A_1$
iff $d(a_1b_1)\geq A_1$. \qed

\blm If $R_1=0$ then condition (ii) of Lemma 2.3 holds for all $N$
with $S_1=0$ iff $\alpha_1\leq 1$.
\elm
\pf By Lemma 2.6, if $S_1=0=R_1$ then $A_1=\alpha_1$ and the condition
(ii) of Lemma 2.3, $d[a_1b_1]\geq A_1$, writes as $d(a_1b_1)\geq
A_1=\alpha_1$.

Since $\ord a_1b_1=R_1+S_1=0$ is even, we have $d(a_1b_1)\geq 1$, so
the condition that $1\geq\alpha_1$ is sufficient. For the necessity,
let $\varepsilon\in\ooo$ with $d(\varepsilon )=1$ and let
$b_1=\varepsilon a_1$. Then $S_1=\ord b_1=\ord a_1=0$ and in
$\fff/\fff^2$ we have $a_1b_1=\varepsilon$ so $d(a_1b_1)=d(\varepsilon
)=1$. Hence we must have $1=d(a_1b_1)\geq\alpha_1$. \qed

\blm (i) If $\alpha_i=0$ then $d[-a_{i,i+1}]\geq 2e$.

(ii) If $\alpha_i=1$ then $-2e<R_{i+1}-R_i\leq 1$ or, equivalently,
either $R_{i+1}-R_i=1$ or $R_{i+1}-R_i$ is even and $2-2e\leq
R_{i+1}-R_i\leq 0$. Moreover, $d[-a_{i,i+1}]\geq R_i-R_{i+1}+1$, with
equality if $R_{i+1}-R_i\neq 2-2e$.

(iii) If $R_{i+1}-R_i\in\{ 2-2e,1\}$ then $\alpha_i=1$
unconditionally. If $2-2e<R_{i+1}-R_i\leq 0$ then $\alpha_i=1$ iff
$d[-a_{i,i+1}]=R_i-R_{i+1}+1$.
\elm
\pf By 1.3, property (2), we have $\alpha_i=0$ iff $R_{i+1}-R_i=-2e$.

Assume that $\alpha_i=1$. By \S1.2, we have $R_{i+1}-R_i\geq -2e$. But
we cannot have $R_{i+1}-R_i=-2e$, since this would imply
$\alpha_i=0$. So $R_{i+1}-R_i>2e$. We cannot have $R_{i+1}-R_i>2e$,
since by \S1.3, property (5), this would imply $\alpha_i>2e$. So
$R_{i+1}-R_i\leq 2e$, which, by property (3), implies $1=\alpha_i\geq
R_{i+1}-R_i$. So $-2e<R_{i+1}-R_i\leq 1$. By \S1.2, $R_{i+1}-R_i$
cannot be odd and negative. So we have either $R_{i+1}-R_i=1$ or
$R_{i+1}-R_i$ is even and $2-2e<R_{i+1}-R_i\leq 0$.

We now use the relation $\alpha_i=\min\{ (R_{i+1}-R_i)/2+e,\,
R_{i+1}-R_i+d[-a_{i,i+1}]\}$. This implies that $\alpha_i\leq
R_{i+1}-R_i+d[-a_{i,i+1}]$, so $d[-a_{i,i+1}]\geq
R_i-R_{i+1}+\alpha_i$, with equality if
$\alpha_i<(R_{i+1}-R_i)/2+e$. If $\alpha_i=0$, so $R_{i+1}-R_i=-2e$,
we get $d[-a_{i,i+1}]\geq R_i-R_{i+1}=2e$, which concludes the proof
of (i). If $\alpha_i=1$ we get $d[-a_{i,i+1}]\geq R_i-R_{i+1}+1$, with
equality if $1<(R_{i+1}-R_i)/2+e$, i.e. if $R_{i+1}-R_i>2-2e$. This
concludes the proof of (ii).

(iii) By \S1.3, properties (3) and (4), if $R_{i+1}-R_i=1$ then
$\alpha_i=R_{i+1}-R_i=1$ and if $R_{i+1}-R_i=2-2e$ then
$\alpha_i=(R_{i+1}-R_i)/2+e=1$. If $2-2e<R_{i+1}-R_i\leq 0$, then the
necessity of $d[-a_{i,i+1}]=R_i-R_{i+1}+1$ follows from
(ii). Conversely, assume that $d[-a_{i,i+1}]=R_i-R_{i+1}+1$. Then
$R_{i+1}-R_i+d[-a_{i,i+1}]=1$ and, since $R_{i+1}-R_i>2-2e$, we also
have $(R_{i+1}-R_i)/2+e>1$. It follows that $\alpha_i=\min\{
(R_{i+1}-R_i)/2+e,\, R_{i+1}-R_i+d[-a_{i,i+1}]\} =1$. \qed

We are interested in the case when $R_1=0$ and $i=1$. We get:

\bco Assume that $m\geq 2$ and $R_1=0$. Then we have:

(i) If $\alpha_1=0$ then $d[-a_{1,2}]\geq 2e$.

(ii) If $\alpha_1=1$ then $-2e<R_2\leq 1$ or, equivalently,
either $R_2=1$ or $R_2$ is even and $2-2e\leq R_2\leq 0$. Moreover,
$d[-a_{1,2}]\geq 1-R_2$, with equality if $R_2\neq 2-2e$.

(iii) If $R_2\in\{ 2-2e,1\}$ then $\alpha_i=1$
unconditionally. If $2-2e<R_2\leq 0$ then $\alpha_1=1$ iff
$d[-a_{1,2}]=1-R_2$.
\eco
 

We define the following statement, which is slightly stronger than
II (a):
\medskip

{\em II (a') $m\geq 3$, $\alpha_1=1$ and $d[-a_{1,2}]=1-R_2$.}
\medskip

(By Corollary 2.9(ii), if $R_1=0$ and $R_2\neq 2-2e$, then the extra
condition that $d[-a_{1,2}]=1-R_2$ is superfluous, as it follows from
$R_1=0$, $\alpha_1=1$.)

\blm Assume that $FM$ is universal and $R_1=0$. Then condition (ii) of Lemma
2.3 holds for every $N$ iff we have I (a) or II (a').
\elm
\pf By Lemma 2.7, the condition (ii) of Lemma 2.3 in the case $S_1=0$
is equivalent to $\alpha_1\leq 1$. We must prove that, assuming that
$R_1=0$ and $\alpha_1\leq 1$, the condition (ii) of Lemma 2.3 holds
for every $N$ with $S_1=1$ iff the additional conditions from II (a'),
$m\geq 3$ and $d[-a_{i,i+1}]=1-R_2$, hold.

Since $\ord a_1b_1=R_1+S_1=1$ is odd, we have $d[a_1b_1]=0$, so
condition (ii) of Lemma 2.3 writes as $0\geq A_1$. If $\alpha_1=0$, so
$R_2=-2e$, then $A_1\leq (R_2-S_1)/2+e=-1/2<0$ so we are done. If
$\alpha_1=1$ and $d[-a_{1,2}]=1-R_2$ then $A_1\leq
R_2-S_1+d[-a_{1,2}]=R_2-1+(1-R_2)=0$ so again we are done.

So we have the sufficiency of the condition $d[-a_{1,2}]=1-R_2$ from II
(a'). For the necessity, assume that $\alpha_1=1$ and $d[-a_{1,2}]\neq
1-R_2$. By Corollary 2.9(ii), this means that $R_2=2-2e$ and
$d[-a_{1,2}]>1-R_2$. Then $R_2-S_1+d[-a_{1,2}]>R_2-1+(1-R_2)=0$ and
$(R_2-S_1)/2+e=((2-2e)-1)/2+e=1/2>0$. It follows that $A_1=\min\{
(R_2-S_1)/2+e,\, R_2-S_1+d[-a_{1,2}]\} >0$, so the condition (ii) of
Lemma 2.3 doesn't hold.

To complete the proof, we show that the remaining condition, $m\geq
3$, from II (a'), follows from the fact that $FM$ is universal. If
$m=2$, then $d(-a_{1,2})=d[-a_{1,2}]=1-R_2<\infty$ so
$-a_{1,2}\notin\fff^2$. It follows that $FM\cong [a_1,a_2]$ is not
isotropic and so it is not universal. \qed

\blm Assume that $R_1=0$ and $R_2=-2e$.

(i) We have $a_{1,2}\in -\fff^2$ or $-\Delta\fff^2$. In the first case
$[a_1,a_2]$ is isotropic. In the second case $[a_1,a_2]$ represents
precisely the elements of $\fff$ with even orders.

In particular, in both cases $[a_1,a_2]$ represents the elements of
$\fff$ with even orders.

(ii) Assume that $m\geq 3$. If $R_3=0$ then $[a_1,a_2,a_3]$ is
isotropic. If $R_3=1$ then $[a_1,a_2,a_3]$ is isotropic iff
$[a_1,a_2]$ is isotropic.
\elm
\pf (i) By 1.2, as a consequence of $R_2-R_1=-2e$, we have
$-a_{1,2}\in\fff^2$ or $\Delta\fff^2$. In the first case $[a_1,a_2]$
is binary of determinant $-1$ so it is isotropic. In the second case,
for every $b\in\fff$ we have $b\rep [a_1,a_2]$ iff
$(a_1b,-a_{1,2})_\p =(a_1b,\Delta )_\p =1$. But this happens iff $\ord
a_1b=\ord b$ is even. (Recall that $\ord a_1=R_1=0$.)

(ii) We have that $[a_1,a_2,a_3]$ is isotropic iff $-a_3\rep
[a_1,a_2]$. If $R_3=0$ then $\ord a_3=R_3$ is even so $-a_3$ is
represented by $[a_1,a_2]$ in both cases from (i). Hence
$[a_1,a_2,a_3]$ is isotropic. Suppose now that $[a_1,a_2,a_3]$ is
isotropic and $R_3=1$. Then $[a_1,a_2]$ represents $-a_3$ and, since
$\ord a_3=R_3$ is odd, this is possible only if $[a_1,a_2]$ is
isotropic. The reverse implication is trivial. \qed

\blm If $m\geq 3$ and $R_1=0$ then $R_3\geq 0$. If moreover $R_2=1$,
then $R_3\geq 1$.
\elm
\pf By the properties of the good BONGs, $R_3\geq R_1=0$. If $R_2=1$
then we cannot have $R_3=0$, since $R_3-R_2$ cannot be odd and
negative. So in this case $R_3\geq 1$. \qed

\blm If $FM$ is universal, $R_1=0$ and $M$ satisfies I (a) or II (a')
then the condition (iii') of Lemma 2.3 is satisfied for every $N$ iff
$M$ satisfies I (b) or II (b), accordingly.
\elm
\pf Suppose that we have I (a). Then $\alpha_1=0$, $R_2=-2e$ and, by
Corollary 2.9(i), $d[-a_{1,2}]\geq 2e$. 

By Lemma 2.11(i), $[a_1,a_2]$ represents all units so if $S_1=\ord
b_1=0$ then $b_1\rep [a_1,a_2]$ so (iii') holds trivially.

Suppose now that $S_1=1$. If $m=2$ then $FM\cong [a_1,a_2]$ must be
isotropic because it is universal. If $R_3\leq S_1=1$ then (iii')
holds trivially. Suppose now that $R_3>1$. Hence $R_3>S_1$ and we also
have $d[-a_{1,2}]+d[-a_{1,3}b_1]\geq d[-a_{1,2}]\geq
2e>2e+S_1-R_3$. Hence condition (iii') holds iff $b_1\rep
[a_1,a_2]$. Since $\ord b_1=S_1$ is odd, by Lemma 2.11(i), this can
only happen if $[a_1,a_2]$ is isotropic. Conversely, if $[a_1,a_2]$ is
isotropic then it is universal so (iii') holds trivially.

In conclusion, the condition (iii') holds iff $M$ satisfies I (b).

Suppose now that we have II (a'). By Lemma 2.12, $R_3\geq 0$ and if
$R_2=1$ then $R_3\geq 1$.

Suppose first that $R_2\leq 0$ and $R_3\leq 1$. We prove that the
condition (iii') of Lemma 2.3 holds unconditionally. By Corollary
2.9(ii), $R_2$ is even and $2-2e\leq R_2\leq 0$. If $R_3=0$ then
$R_3\leq S_1$ so (iii') holds trivially. Suppose now that $R_3=1$. If
$S_1=1$ then $R_3\leq S_1$ so (iii') holds tivially. If $S_1=0$ then
$\ord a_{1,3}b_1=R_1+R_2+R_3+S_1$ is odd so $d[-a_{1,3}b_1]=0$. (We
have $R_1=S_1=0$, $R_3=1$ and $R_2$ is even.) Hence
$d[-a_{1,2}]+d[-a_{1,3}b_1]=1-R_2+0\leq  2e-1=2e+S_1-R_3$ so again
(iii') holds trivially.

So we are left with the case when $R_2=1$ or $R_3>1$.

We claim that $d(-a_{1,2})=d[-a_{1,2}]=1-R_2$. If $R_2=1$ then $\ord
a_{1,2}=R_1+R_2=1$ so $d(-a_{1,2})=0=1-R_2$. Suppose now that $R_2\leq
0$ and $R_3>1$. Since $R_2\geq 2-2e$, we have $1-R_2\leq 2e-1$. If
$R_3-R_2>2e$ then, by 1.3, property (5), we have
$\alpha_2>2e>1-R_2$. If $R_3-R_2\leq 2e$ then, by the property (3),
we have $\alpha_2\geq R_3-R_2>1-R_2$. So we have $\min\{
d(-a_{1,2}),\alpha_2\} =d[-a_{1,2}]=1-R_2$ and $\alpha_2>1-R_2$. It
follows that $d(-a_{1,2})=1-R_2$.

Suppose first that $S_1\equiv R_2+R_3\pmod 2$. Recall that if $R_2=1$
then $R_3\geq 1$ and if $R_2\leq 0$ then $R_3>1$. So, with the
exception of the case $R_2=R_3=1$, we have $R_3>1\geq S_1$. If
$R_2=R_3=1$ then $S_1\equiv R_2+R_3\equiv 0\pmod 2$ and $S_1\in\{
0,1\}$ so $S_1=0<R_3$. Hence in this case the inegality $R_3>S_1$ from
Lemma 2.3(iii') is satisfied. Condition (iii') states that if moreover
$d[-a_{1,2}]+d[-a_{1,3}b_1]>2e+S_1-R_3$, then $b_1\rep  [a_1,a_2]$.

Note that $d[-a_{1,3}b_1]=\min\{ d(-a_{1,3}b_1),\alpha_3\}$, with
$\alpha_3$ ignored if $m=3$.

If $m\geq 4$ then $d[-a_{1,3}b_1]\leq\alpha_3$ so if
$d[-a_{1,2}]+\alpha_3\leq 2e+S_1-R_3$ then also
$d[-a_{1,2}]+d[-a_{1,3}b_1]\leq 2e+S_1-R_3$ so condition (iii') holds
trivially.

Suppose now that $m=3$ or $m\geq 4$ and
$d[-a_{1,2}]+\alpha_3>2e+S_1-R_3$. Then, by the formula
$d[-a_{1,3}b_1]=\min\{ d(-a_{1,3}b_1),\alpha_3\}$,
$d[-a_{1,2}]+d[-a_{1,3}b_1]>2e+S_1-R_3$ is equivalent to
$d[-a_{1,2}]+d(-a_{1,3}b_1)>2e+S_1-R_3$. So, for (iii') to hold,
$[a_1,a_2]$ must represent every $b_1\in\fff$ with $\ord b_1=S_1$ and
$d[-a_{1,2}]+d(-a_{1,3}b_1)>2e+S_1-R_3$. We have $R_2\geq 2-2e$ so
$d(-a_{1,2})=1-R_2\leq 2e-1$. Then, by \S1.1, that there is
$\varepsilon\in\ooo$ such that $d(\varepsilon
)=2e-d(-a_{1,2})=2e-d[-a_{1,2}]$ and
$(\varepsilon,-a_{1,2})_\p =-1$. Since $\ord
a_{1,3}=R_1+R_2+R_3=R_2+R_3\equiv S_1\pmod 2$, there is $b\in\fff$,
say, $b=-\pi^{S_1-R_2-R_3}a_{1,3}$, such that $\ord b=S_1$ and $b\in
-a_{1,3}\fff^2$. It follows that $-a_{1,3}b\in\fff^2$ and
$-a_{1,3}b\varepsilon\in\varepsilon\fff^2$ so
$d(-a_{1,3}b)=\infty$ and $d(-a_{1,3}b\varepsilon
)=d(\varepsilon )=2e-d[-a_{1,2}]$. In both cases when $b_1=b$ or
$b\varepsilon$, we have $\ord b_1=\ord b=S_1$ and $d(-a_{1,3}b_1)\geq
2e-d[-a_{1,2}]$, which implies $d[-a_{1,2}]+d(-a_{1,3}b_1)\geq
2e>2e+S_1-R_3$. So in both cases we have ${b_1\rep [a_1,a_2]}$,
which is equivalent to $(a_1b_1,-a_{1,2})_\p =1$. We get
$(a_1b,-a_{1,2})_\p =(a_1b\varepsilon,-a_{1,2})_\p =1$, which implies
$(\varepsilon,-a_{1,2})_\p =1$. But this contradicts the choice of
$\varepsilon$. 

So the condition that $m\geq 4$ and $d[-a_{1,2}]+\alpha_3\leq
2e+S_1-R_3$ is not only sufficient, but also necessary. Since
$d[-a_{1,2}]=1-R_2$, this inequality writes as $\alpha_3\leq
2e-1+S_1+R_2-R_3$. We have $S_1\equiv R_3-R_2\pmod 2$ and $S_1\in\{
0,1\}$, so $S_1=R_3-R_2-2[\frac{R_3-R_2}2]$. It follows that
$2e-1+S_1+R_2-R_3=2(e-[\frac{R_3-R_2}2])-1$. 

So we have proved that if $R_2=1$ or $R_3>1$ then the condition that
$m\geq 4$ and $\alpha_3\leq 2(e-[\frac{R_3-R_2}2])-1$ is necessary
and sufficient for condition (iii') of Lemma 2.3 to hold in the case
when $S_1\equiv R_2+R_3\pmod 2$. To conclude the proof, we show that
this condition is also sufficient for (iii') to hold in the case when
$S_1\equiv R_2+R_3+1\pmod 2$. So assume that $\alpha_3\leq
2(e-[\frac{R_3-R_2}2])-1$ and $S_1\equiv R_2+R_3+1\pmod 2$. Suppose
that $d[-a_{1,2}]+d[-a_{1,3}b_1]>2e+S_1-R_3$. We have
$d[-a_{1,2}]=1-R_2$ and $\ord
a_{1,3}b_1=R_1+R_2+R_3+S_1=R_2+R_3+S_1\equiv 1\pmod 2$ so
$d[-a_{1,3}b_1]=0$. Hence $(1-R_2)+0>2e+S_1-R_3\geq 2e-R_3$, which
implies $R_3-R_2>2e-1$, so $R_3-R_2\geq 2e$. It follows that
$\alpha_3\leq 2(e-[\frac{R_3-R_2}2])-1\leq
2(2-[\frac{2e}2])-1=-1$. But $\alpha_3\geq 0$, by \S1.3, property
(2). Contradiction. Hence $d[-a_{1,2}]+d[-a_{1,3}b_1]\leq 2e+S_1-R_3$
and so (iii') holds trivially. \qed

\blm If $FM$ is universal, $R_1=0$ and $M$ satisfies I (a) or II (a')
then the condition (iv) of Lemma 2.3 is satisfied for every $N$ iff
$M$ satisfies I (c) or II (c), accordingly.
\elm
\pf Recall that $b_1\rep [a_1,a_2,a_3]$ iff $b_1\notin
-a_{1,3}\fff^2$ or $[a_1,a_2,a_3]$ is isotropic.

Let $S\in\{ 0,1\}$ such that $S\equiv R_2+R_3\pmod 2$. We claim that
condition (iv) of Lemma 2.3 holds for every $N$ iff the following
statement holds.

(*) If $m\geq 3$, $R_3\leq S$ and either $m=3$ or $R_4-S>2e$ then
$[a_1,a_2,a_3]$ is isotropic.

Assume first that $m\geq 4$.

If $S_1\equiv R_2+R_3+1\pmod 2$ then $\ord b_1=S_1$ and $\ord
a_{1,3}=R_1+R_2+R_3=R_2+R_3$ have opposite parities so we cannot have
$b_1\in -a_{1,3}\fff^2$. Therefore $b_1\rep [a_1,a_2,a_3]$ so in
this case condition (iv) of Lemma 2.3, holds unconditionally.

Suppose now that $S_1\equiv R_2+R_3\pmod 2$, i.e. that $S_1=S$. If
$R_3>S_1=S$ or $R_4-S=R_4-S_1\leq 2e$, then (iv) holds trivially. So
we assume that $R_3\leq S=S_1$ and $R_4-S=R_4-S_1>2e$. Then
condition (iv) of Lemma 2.3 states that $b_1\rep [a_1,a_2,a_3]$. So
$[a_1,a_2,a_3]$ must represent all elements of $\fff$ of order
$S_1=S$. Since $\ord a_{1,3}=R_1+R_2+R_3=R_2+R_3\equiv S\pmod 2$ there
is $b_1\in -a_{1,3}\fff^2$ with $\ord b_1=S$, say,
$b_1=-\pi^{S-R_2-R_3}a_{1,3}$. Then $b_1\rep [a_1,a_2,a_3]$ and
$b_1\in -a_{1,3}\fff^2$, so $[a_1,a_2,a_3]$ is isotropic.

Conversely, if $[a_1,a_2,a_3]$ is isotropic, then it is universal, so
(iv) holds. 

If $m\leq 3$ then (iv) is vacuous, but the case when $m=3$ and
$R_3\leq S$ can still be included here because when $m=3$ we have that
$FM\cong [a_1,a_2,a_3]$ is universal, so isotropic.

Suppose first that we have I (a). If $R_3>1$ then $R_3>S$ so (*) holds
trivially. If $R_3=0$ then, by Lemma 2.11(ii), $[a_1,a_2,a_3]$ is
isotropic, so (*) holds unconditionally. We are left with the case
$R_3=1$. Then $R_2+R_3=1-2e$ is odd so $S=1$ and $R_3\leq S$
holds. The inequality $R_4-S>2e$ writes as $R_4>2e+1$. Therefore (*)
is equivalent to:

I (c') If $m\geq 3$, $R_3=1$ and either $m=3$ or $R_4>2e+1$, then
$[a_1,a_2,a_3]$ is isotropic. 

But when $R_3=1$, by Lemma 2.11(ii), $[a_1,a_2,a_3]$ is isotropic iff
$[a_1,a_2]$ is isotropic. Hence I (c') is equivalent to I (c).

Assume now that $M$ satisfies II (a').

Suppose first that $R_3>1$. Since $S\leq 1$ we get $R_3>S$ so (*)
holds trivially. So we may assume that $R_3\leq 1$. Next suppose that
$R_2=1$. By Lemma 2.12, $R_3\geq 1$, so $R_3=1$. Since $R_2+R_3=2$ is
even, we get $S=0$ and again $R_3>S$, so (*) holds trivially.

Since for $R_2=1$ or $R_3>1$ (*) holds unconditionally, we are left
with the case when $R_2\leq 0$ and $R_3\leq 1$. Since $R_2\leq 0$ we
have that $R_2$ is even and so $S\equiv R_2+R_3\equiv R_3\pmod
2$. Since $R_3,S\in\{ 0,1\}$ and $R_3\equiv S\pmod 2$, we have
$S=R_3$. In particular, $R_3\geq S$ holds. So, in order that (*) holds
we need that $[a_1,a_2,a_3]$ is isotropic if $m=3$ or $m\geq 4$ and
$R_4-R_3=R_4-S>2e$. So (*) writes as follows. If $m\geq 3$, $R_2\leq
0$ and $R_3\leq 1$ and either $m=3$ or $m\geq 4$ and $R_4-R_3>2e$,
then $[a_1,a_2,a_3]$ is isotropic. But the condition that $m\geq 3$ is
part of II (a') so it can be dismissed. So we get II (c). \qed

{\bf Proof of Theorem 2.1.} By Lemmas 2.2, 2.5, 2.10, 2.13 and 2.14,
$M$ is universal iff $FM$ is universal, $R_1=0$ and we have either I
(a), (b) and (c) or II (a'), (b) and (c). Since $FM$ is universal, we
have $m\geq 2$. Then, to conclude the proof, we must show that II (a')
can be replaced by II (a) and that the condition that $FM$ is
universal is superfluous, as if $m\geq 2$ and $R_1=0$ then it follows
both from I and from II.

First we prove that the extra condition from II (a'), that
$d[-a_{1,2}]=1-R_2$, is superfluous, as it follows from II (a) and
(b). Since $R_1=0$ and $\alpha_1=1$, by Corollary 2.(ii),
$d[-a_{1,2}]\geq 1-R_2$, with equality if $R_2\neq 2-2e$. So we only
have to consider the case $R_2=2-2e$, when we only have
$d[-a_{1,2}]\geq 1-R_2=2e-1$. Assume that $R_3>1$. Then, by II (b), we
have that $m\geq 4$ and $\alpha_3\leq 2(e-[\frac{R_3-R_2}2])-1$. But
$R_3\geq 2$ so $R_3-R_2=R_3-(2-2e)\geq 2e$. Then $\alpha_3\leq
2(e-[\frac{R_3-R_2}2])-1\leq 2(e-[\frac{2e}2])-1=-1$. But, by \S1.3,
property (2), we have $\alpha_3\geq 0$. Contradiction. Hence $R_3\leq
1$. If $R_3=0$ then $R_3-R_2=2e-2$ so, by \S1.3, property (4),
$\alpha_2=(R_3-R_2)/2+e=2e-1$. If $R_3=1$ then $R_3-R_2=2e-1$ is odd
and $<2e$ so, by the property (3), $\alpha_2=R_3-R_2=2e-1$. Then
$2e-1\leq d[-a_{1,2}]=\min\{ d(-a_{1,2}),\alpha_2\}\leq\alpha_2=2e-1$
so $d[-a_{1,2}]=2e-1=1-R_2$.

Next we prove that if $m\geq 2$, $R_1=0$ and we have I or II, then
$FM$ is universal. If $m=2$ then we are in case of I and we have
by I (b) that $FM=[a_1,a_2]$ is isotropic and so it is
universal. Suppose now that $m=3$. We consider first case I. If
$R_3>1$ then $[a_1,a_2,a_3]$ is isotropic by I (b), if $R_3=1$ then it
is isotropic by I (c) and if $R_3=0$ it is isotropic by Lemma
2.11(ii). If we are in the case II then we cannot have $R_2=1$ or
$R_3>1$, since by II (b) this implies that $m\geq 4$. In the remaining
case, $R_2\leq 0$ and $R_3\leq 1$, since $m=3$ we have by I (c) that
$[a_1,a_2,a_3]$ is isotropic. So, in all cases, $FM=[a_1,a_2,a_3]$ is
isotropic and so universal. If $m\geq 4$ then $FM$ is universal
unconditionally. \qed

{\bf Remark.} If we make the convention that $R_i\gg 0$ for $i>m$ then
conditions I and II can be written in a more compact way, without
refference to the value of $n$. Namely one can write them as:

{\em I (a) $\alpha_1=0$ or, equivalently, $R_2=-2e$.

~~(b) If $R_3>1$, then $[a_1,a_2]$ is isotropic.

~~(c) If $R_3=1$ and $R_4>2e+1$, then $[a_1,a_2]$ is isotropic.}

\medskip

{\em II (a) $\alpha_1=1$.

~~(b) If $R_2=1$ or $R_3>1$, then $m\geq 4$ and $\alpha_3\leq
2(e-[\frac{R_3-R_2}2])-1$.

~~(c) If $R_2\leq 0$, $R_3\leq 1$ and $R_4-R_3>2e$, then
$[a_1,a_2,a_3]$ is isotropic.}

\section{Main result in terms of Jordan decompostions}

We now give, without a proof, a translation of Theorem 2.1, in terms
of Jordan decompositions.

\btm Let $M=M_1\perp\cdots\perp M_t$ be a Jordan decompostion, with
$\sss M_k=\p^{r_k}$, $\nnn M^{\sss M_k}=u_k$ and
$\www M^{\sss M_k}=\www_k$ for $1\leq k\leq t$. For $1\leq k\leq t-1$ we
consider
the ideal $\FFF_k$ defined in [OM, \S93E.]. Then $M$ is universal if
$\rank M\geq 2$, $\nnn M=\oo$, or, equivalently, $u_1=0$, and one of
the following happens:

(1) $\rank M_1\geq 4$ and $\www_1\supseteq\p$.

(2) $\rank M_1=3$, $\www_1=\p$ and one of the following happens:

(2.1) $t\geq 2$ and $u_2\leq 2e$.

(2.2) $M_1$ is isotropic.

(3) $\rank M_1=2$ and one of the following happens:

(3.1) $\sss M_1=\frac 12\oo$ and one of the following happens:

(3.1.1) $t\geq 2$ and $u_2=0$.

(3.1.2) $t\geq 2$, $u_2=1$ and either $\rank M_2\geq 2$ or $\rank
M_2=1$, $t\geq 3$ and $u_3\leq 2e+1$.

(3.1.3) $M_1\cong\frac 12 A(0,0)$.

(3.2) $\www M=\p$, $m\geq 3$ and one of the following happens:

(3.2.1) $u_2>1$, $\rank M_2\geq 2$ and $\www_2\supset
4\p^{r_1+u_2-2[u_2/2]}$.

(3.2.2) $u_2>1$, $\rank M_2=1$, $t\geq 3$ and $\FFF_2\supset
4\p^{r_1-2[u_2/2]}$.

(3.2.3) $u_2\leq 1$ and $\rank M_2\geq 2$.

(3.2.4) $u_2\leq 1$, $\rank M_2=1$, $t\geq 3$ and $u_3\leq u_2+2e$.

(3.2.5) $u_2\leq 1$, $\rank M_2=1$ and $M_1\perp M_2$ is isotropic.

(4) $\rank M_1=1$, $u_2=1$, $m\geq 4$ and one of the following
happens:

(4.1) $\rank M_2\geq 3$.

(4.2) $\rank M_2=2$ and $u_3\leq 2e$.

(4.3) $\rank M_2=1$ and one of the following happens:

(4.3.1) $\rank M_3\geq 2$ and $\www_3\supset 4\p^{u_3-2[(u_3-1)/2]}$.

(4.3.2) $\rank M_3=1$ and $\FFF_3\supset 4\p^{-2[(u_3-1)/2]}$.
\etm

\section{Some remarks on $n$-universality}

Here we consider the more general problem of $n$-universality. An
integral lattice $M$ is said to be $n$-universal if it represents all
integral lattices of rank $n$. We will not solve this problem, but we
show how it can be reduced, by an inductive argument, to the case when
$n\leq 4$.

We start by noting that every integral lattice $N$ is included in an
$\oo$-maximal lattice over $FN$. (See [OM, 82:18].) Hence it suffices
to consider only the $\oo$-maxial lattices, i.e. we have:

\blm A lattice $M$ is $n$-universal iff it represents all
$\oo$-maximal lattices.
\elm

Incidentally, this proves that the number of lattices $N$ to be
considered is finite. Indeed, by [OM, Theorem 91.2], for any $n$-ary
quadratic space $W$ we have, up to an isometry, only one $\oo$-maximal
lattice over $W$. Hence the number of classes of $n$-ary maximal
lattices is equal to the number of classes of $n$-ary quadratic
spaces, which is finite.

For convenience, we denote by $H$ a hyperbolic plane.

By [OM, 82:23], if $N$ is maximal over an isotropic quadratic space
then $N=J\perp N'$, where $FJ\cong H$. Since $N$ is $\oo$-maximal, so
are $J$ and $N'$. By [OM, 93:11], we have $J\cong\frac 12A(0,0)$. By
induction, we get

\blm If $N,N'$ are $oo$-maximal lattices over the quadratic spaces
$W,W'$ and $W\cong H^k\perp W'$ for some $k\geq 0$, then
$N\cong\frac 12A(0,0)^k\perp N'$.

(Here by $H^k$ we mean an orthogonal sum of $k$ copies of
$H$. Similarly for $\frac 12A(0,0)^k$.)
\elm

This allows us to determine $\oo$-maximal lattices over isotropic
vector spaces from $\oo$-maximal lattices over spaces of smaller
dimensions.

\blm Let $k\geq 1$ and let $M$ and $N'_i$, with $i\in I$, be integral lattices.
Then $M$ represents $\frac 12 A(0,0)^k\perp N'_i$ $\forall i\in I$ if and only
if $M\cong\frac 12 A(0,0)^k\perp M'$ for some integral lattice $M'$ that
represents $N'_i$ $\forall i\in I$.

Here $I$ is a non-empty set of indices.
\elm
\pf The ''if'' part of the statemen is trivial. We now prove the
''only if'' part.

For every $i\in I$ we have $\frac 12A(0,0)^k\perp N'_i\rep M$, so $M$ contains
a sublattice $J_i\perp N''_i$, where $J_i\cong\frac 12A(0,0)^k$ and
$N''_i\cong N'_i$. Since $M$ is integral, we have a $\nnn M\subseteq\oo$ and so
$\sss M\subseteq\frac 12\oo$. Since $M$ contains the $\frac 12\oo$-modular
lattice $J_i$, we must have $\sss M=\frac 12\oo$ and $M$ splits $J_i$, i.e.
$M=J_i\perp M'_i$ for some lattice $M'_i$. Then
$J_i\perp N''_i\subseteq M=J_i\perp M'_i$ implies that $N''_i\subseteq M'_i$.

For every $i,j\in I$ we have $M=J_i\perp M'_i=J_j\perp M'_j$. Since
$J_i\cong J_j\cong\frac 12A(0,0)^k$, by the cancelattion law
[OM, 93:14], we have $M'_i\cong M'_j$. Let $M'$ be the common class of
all $M'_i$ with $i\in I$. Then for every $i\in I$ we have
$M=J_i\perp M'_i\cong\frac 12 A(0,0)^k\perp M'$. Also
$N'_i\cong N''_i\subseteq M'_i\cong M'$, so $N'_i\rep M'$. \qed

We denote by $m(V)$ the Witt index of a quadratic space $V$, which is the
dimension of the largest isotropic subspace of $V$ or, equivalently, $m(V)$ is
maximal with the property that $V$ splits $H^{m(V)}$.

\blm Let $k\geq 1$ and let $n\geq 2k+1$. If $M$ is an integral lattice then the
following are equivalent:

(i)  $M$ represents all integral $n$-ary lattices over spaces with
Witt index $\geq k$.

(ii) $M\cong\frac 12A(0,0)^k\perp M'$ for some $n-2k$-universal lattice $M'$.
\elm
\pf Statement (i) is equivalent to 

{\it (i') $M$ represents all $\oo$-maximal $n$-ary lattices over
spaces with Witt index $\geq k$.}

The set of all $n$-ary quadratic forms $W$ with $m(W)\geq k$ can be
written as
$$\{ H^k\perp W',\, W'\text{ is an }n-2k\text{-ary quadratic form}\}.$$

Then by Lemma 4.2, the set of all $\oo$-maximal $n$-ary lattices over
spaces with Witt index $\geq k$ can be written as
$$\{\frac 12A(0,0)^k\perp N',\, N'\text{ is an }\oo\text{-maximal }
n-2k\text{-ary quadratic lattice}\}.$$

Then if in Lemma 4.3 we take $\{ N'_i\mid i\in I\}$ to be the set of
all $n-2k$-ary $\oo$-maximal lattices, the condition (i') holds iff
$M=\frac 12A(0,0)^k\perp M'$, where $M'$ represents all $n-2k$-ary
$\oo$-maximal lattices, which, by Lemma 4.1, is equivalent to $M'$
being $n-2k$-universal. This concludes the proof. \qed

\bco (i) $M$ is $3$-universal iff $M=\frac 12A(0,0)\perp M'$, where
$M'$ is an $1$-universal lattice, and $M$ represents all ternary
$\oo$-maximal lattices over anisotropic quadratic spaces.

(ii) $M$ is $4$-universal iff $M=\frac 12A(0,0)\perp M'$, where $M'$
is a $2$-universal lattice, and $M$ represents the $\oo$-maximal
lattice over $[1,-\Delta,\pi,-\Delta\pi ]$.

(iii) If $k\geq 1$, then $M$ is $2k+3$-universal iff
$M\cong \frac 12A(0,0)^k\perp M'$, where $M'$ is a $3$-universal
lattice.

(iv) If $k\geq 1$, then $M$ is $2k+4$-universal iff
$M\cong \frac 12A(0,0)^k\perp M'$, where $M'$ is a $4$-universal
lattice.
\eco
\pf (i) By Lemma 4.4, the condition that $M=\frac 12A(0,0)\perp M'$,
where $M'$ is an $1$-universal lattice, is equivalent to $M$
representing all integral lattices over ternary quadratic spaces $W$
with $m(W)\geq 1$, i.e. with $W$ isotropic. In order that $M$ is
$3$-universal, it must also represent the integral lattices over
ternary anisotropic quadratic spaces or, equvalently, all maximal
$\oo$-lattices over ternary anisotropic quadratic spaces.

(ii) Similarly the condition that $M=\frac 12A(0,0)\perp M'$, where
$M'$ is a $2$-universal lattice is equivalent to $M$ representing all
integral quadratic forms over isotropic quaternary spaces. What
remains is the only anisotropic quaternary quadratic space,
i.e. $W=[1,-\Delta,\pi,-\pi\Delta ]$.

Finally, if $\dim W\geq 5$, then $W$ is isotropic. So every lattice
$W$ writes as $W\cong H^k\perp W'$, with $\dim W'\leq 4$. So if
$\dim W=2k+3$ or $2k+4$, then $W$ splits at least $k$ copies of $H$
and so $m(W)\geq k$. 

Then, by Lemma 4.4, the condition that
$M\cong \frac 12A(0,0)^k\perp M'$, where $M'$ is a $n-2k$-universal 
lattice, from (iii) and (iv) is equivalent to $M$ representing all
integral lattices over $n$-ary quadratic spaces $W$ with
$m(W)\geq k$. But if $n=2k+3$ or $2k+4$ then $m(W)\geq k$ so all
quadratic spaces are covered. Hence our condition is equivalent to $M$
being $n$-universal. \qed

Corollary 4.5 reduces the problem to determining all $n$-maximal
lattices for $n=1,2,3,4$. It also gives some necessary conditions in
the case $n=3$ or $4$, derived from those for $n=1$ and $2$,
respectively. However translating conditions these conditions in terms
of good BONGs is not that straightforward.

\blm A $2k$-dimesional lattice $J$ with $R_i(J)=R_i$ is
$\frac 12\oo$-modular with norm $\oo$ if and only if the sequence
$R_1,\ldots,R_{2k}$ is $0,-2e,\ldots,0,-2e$. 

Moreover we have $\det FJ=(-1)^k$ if $FJ\cong\frac 12A(0,0)^k$ and
$\det FJ=(-1)^k\Delta$ if
$\det FJ=\frac 12A(0,0)^{k-1}\perp\frac 12A(2,2\rho )$. 
\elm
\pf We use [B3, Lemma 2.13]. (See also the proof of
[B1, Lemma 4.7].) Since $J$ has only one Jordan component of rank $2k$
and $\sss J=\p^r$, $\nnn J=\p^u$, where $r=-e$ and $u=0$, the sequence  
$R_1,\ldots,R_{2k}$ is $u,2r-u,\ldots,u,2r-u$,
i.e. $0,-2e,\ldots,0,-2e$.

The second statement follows from $\det FA(0,0)=\det H=-1$ and
$\det FA(2,2\rho )=\det [1,-\Delta ]=-\Delta$. \qed

\blm If $J\cong\prec a_1,\ldots,a_{2k}\succ$,
$M'\cong\prec a_{2k+1},\ldots,a_m\succ$ relative to good BONGs, $J$ is
$\frac 12\oo$-modular of norm $\oo$ and $M'$ is integral, then
$J\perp M'\cong\prec a_1,\ldots,a_m\succ$ relative to a good BONG.
\elm
\pf Let $R_i=\ord a_i$. Since $M'$ is integral, we have
$R_{2k+1}=R_1(M')\geq 0$. If $m\geq 2k+2$, then also
$R_{2k+2}-R_{2k+1}\geq -2e$ and so
$R_{2k+2}\geq R_{2k+1}-2e\geq -2e$. Then $R_{2k}=-2e<0\leq R_{2k+1}$,
$R_{2k-1}=0\leq R_{2k+1}$ and, if $m\geq 2k+2$,
$R_{2k}=-2e\leq R_{2k+2}$. Then our statement follows from 
[B1, Corollary 4.4(v)]. \qed 

\blm An integral lattice $M\cong\prec a_1,\ldots,a_m\succ$ relative to
a good BONG, with $R_i(M)=R_i$ splits $\frac 12A(0,0)^k$ if and only
if the following hold:

(i) $R_1=R_3=\cdots =R_{2k-1}=0$ and $R_2=R_4=\cdots =R_{2k}=-2e$.

(ii) Either $d[(-1)^ka_{1,2k}]>2e$ or $m\geq 2k+1$ and $R_{2k+1}=0$.
\elm
\pf For the ''only if'' part, we note that the $R_i$'s and
$d[(-1)^ka_{1,2k}]$ are invariants independent of the choice of the
good BONG. So we may assume that $a_1,\ldots,a_m$ are those from Lemma
4.7. Then $J=\prec a_1,\ldots,a_{2k}\succ\cong\frac 12A(0,0)^k$, which
implies that the sequence $R_1,\ldots,R_{2k}$ is $0,-2e,\ldots,0,-2e$ and
$a_{1,2k}=\det FJ=(-1)^k$ in $\fff/\ffs$. In particular, we have
(i). Also in $\fff/\ffs$ we have $(-1)^ka_{1,2k}=1$ so
$d((-1)^ka_{1,2k})=\infty$. If $m=2k$ then
$d[(-1)^ka_{1,2k}]=d((-1)^ka_{1,2k})=\infty$ so (ii) holds. If
$m\geq 2k+1$ then $R_{2k+1}\geq R_{2k-1}=0$. If $R_{2k+1}=0$ then (ii)
holds. If $R_{2k+1}>0$ then $R_{2k+1}-R_{2k}>0-(-2e)=2e$ so
$\alpha_{2k}>0$. Together with $d((-1)^ka_{1,2k})=\infty$, this
implies that
$d[(-1)^ka_{1,2k}]=\min\{ d((-1)^ka_{1,2k}),\alpha_{2k}\} >2e$, so
(ii) holds.

For the ''if'' part, the condition (i) implies, by Lemma 4.6, that if
$J=\prec a_1,\ldots,a_{2k}\succ$, then $J$ is $\frac 12\oo$-modular of
norm $\oo$. Moreover, $a_{1,2k}=\det J=(-1)^k$ or $(-1)^k\Delta$,
corresponding to $J\cong\frac 12A(0,0)^k$ or
$\frac 12A(0,0)^{k-1}\perp\frac 12A(2,2\rho )$, respectively. If
$m\geq 2k+1$ then $R_{2k+1}\geq R_{2k-1}=0>-2e=R_{2k}$ so, by
[B1, Corollary 4.4(i)], we have $M=J\perp M'$, where
$M'=\prec a_{2k+1},\ldots,a_m\succ$. We now use property (ii). If
$d[(-1)^ka_{1,2k}]>2e$ then $d((-1)^ka_{1,2k})>2e$ so in $\fff/\ffs$
we have $(-1)^ka_{1,2k}=1$, i.e. $a_{1,2k}=(-1)^k$, and so
$J\cong\frac 12A(0,0)^k$, which concludes the proof. Suppose now that
$R_{2k+1}=0$. If $a_{1,2k}=(-1)^k$ in $\fff/\ffs$ then again
$J\cong\frac 12A(0,0)^k$ and we are done. So assume that
$a_{1,2k}=(-1)^{2k}\Delta$. Now $\ord a_{2k}a_{2k+1}=-2e$ is even and
so $d(-a_{2k}a_{2k+1})>0=e-(R_{2k+1}-R_{2k})/2$. Then, by
[B1, Definition 6(iii)],
$g(a_{2k+1}/a_{2k})=(1+\p^{(R_{2k+1}-R_{2k})/2+e})\oos
=(1+\p^{2e})\oos=\langle\Delta\rangle\oos$. Since
$\Delta\in g(a_{2k+1}/a_{2k})$, by [B1, 3.12 and Lemma 4.9(ii)], we
have $\prec a_{2k},a_{2k+1}\succ\cong\prec\Delta a_{2k},\Delta
a_{2k+1}\succ$ and $M\cong\prec a_1,\ldots,\Delta a_{2k},\Delta
a_{2k+1},\ldots,a_m\succ$ relative to some good BONG. Then we replace
the splitting $M=J\perp M'$ by  $M=J'\perp M''$, where
$J'=\prec a_1,\ldots,\Delta a_m\succ$ and 
$M''=\prec\Delta a_{2k+1},\ldots,a_m\succ$. Since
$\det J=(-1)^k\Delta$, we get $\det J'=\Delta\det J=(-1)^k$, so this
time $J'\cong\frac 12A(0,0)^k$. \qed 

\blm Let $M,M'$ be integral lattices with
$M\cong\prec a_1,\ldots,a_m\succ$ and 
$M'\cong\prec a'_1,\ldots,a'_{m'}\succ$ relative to some good BONGs
such that $M\cong\frac 12A(0,0)^k\perp M'$ for some $k\geq 1$. Let
$R_i(M)=R_i$, $\alpha_i(M)=\alpha_i$, $R_i(M')=R'_i$ and
$\alpha_i(M')=\alpha'_i$, Then we have:

(i) $m'=m-2k$, $R'_i=R_{2k+i}$ for $1\leq i\leq m'$ and
$\alpha'_i=\alpha_{2k+i}$ for $1\leq i\leq m'-1$.

(ii) If $\e\in\fff$ then $d[\e a'_{i,j}]=d[\e a_{2k+i,2k+j}]$ for every
$1<i\leq j\leq m'$ and $d[\e a'_{1,j}]=d[(-1)^k\e a_{1,2k+j}]$ for
every $1<j\leq m'$.

(iii) The BONG of $M'$ can be chosen with the property that
$[a_1,\ldots,a_{2k+i}]\cong H^k\perp [a'_1,\ldots,a'_i]$ for
$1\leq i\leq m'$ 
\elm
\pf The $m'=m-2k$ part of (i) follows by considering dimesions in
$M\cong\frac 12A(0,0)^k\perp M'$. 

Since $R_i$, $\alpha_i$ and $d[\e a_{i,j}]$ are invariants of the
lattice $M$, the statements (i) and (ii) are independent of the choice
of the good  BONG of $M$. So we may assume that $a_1,\ldots,a_m$ are
obtained as in Lemma 4.7. Namely, we take $a_1,\ldots,a_{2k}$ such
that $\prec a_1,\ldots,a_{2k}\succ\cong\frac 12A(0,0)^k$ and we take
$a_{2k+i}=a'_i$ for $1\leq i\leq m'$, so that
$M'\cong\prec a'_1,\ldots,a'_{m'}\succ =\prec a_{2k+1},\ldots,a_m\succ$.
By considering orders in $a'_i=a_{2k+i}$, we get the $R'_i=R_{2k+i}$
relation. We have $[a_1,\ldots,a_{2k}]\cong H^k$ so in $\fff/\ffs$ we
have $a_{1,2k}=(-1)^k$. From this we conclude that for
$1\leq i\leq m'$ in $\fff/\ffs$ we have
$a_{1,2k+i}=a_{1,2k}a_{2k+1,2k+i}=(-1)^ka'_{1,i}$.

For the relation $\alpha'_i=\alpha_{2k+i}$ we use the following
formula:
$$\alpha_{2k+i}=\min\{\alpha'_i,R_{2k+i+1}-R_{2k+1}+\alpha_{2k}\}.$$
This happens because, by [B3, Lemma 2.1],
$\alpha'_i=\alpha_i(\prec a_{2k+1},\ldots,a_m\succ )$ can replace
$(R_{2k+i+1}-R_{2k+i})/2+e$ and all the terms in the definition of
$\alpha_{2k+i}$ corresponding to indices $2k+1\leq j\leq m-1$ and, by
[B3, Lemma 2.4(i)], $R_{2k+i+1}-R_{2k+1}+\alpha_{2k}$ can replace all
terms corresponding to $1\leq j\leq 2k$. But
$R_{2k+1}-R_{2k}\geq R_{2k-1}-R_{2k}=0-(-2e)=2e$ so, by
[B3, Lemma 2.7(ii)], $\alpha_{2k}=(R_{2k+1}-R_{2k})/2+e$. It follows
that
$R_{2k+i+1}-R_{2k+1}+\alpha_{2k}=R_{2k+i+1}-(R_{2k}+R_{2k+1})/2+e$. But
the sequence $(R_j+R_{j+1})$ is increasing, so 
$R_{2k}+R_{2k+1}\leq R_{2k+i}+R_{2k+i+1}$. Then
\begin{multline*}
R_{2k+i+1}-R_{2k+1}+\alpha_{2k}\geq
R_{2k+i+1}-(R_{2k+i}+R_{2k+i+1})/2+e \\ 
=(R_{2k+i+1}-R_{2k+i})/2+e=(R'_{i+1}-R'_i)/2+e\geq\alpha'_i.
\end{multline*}
Hence $\alpha_{2k+i}=\alpha'_i$, which concludes the proof of (i).

(ii) If $1\leq j\leq m'-1$ then
$d[\e a'_{1,j}]=\min\{ d(\e a'_{1,j}),\alpha'_j\}$ and
$d[(-1)^k\e a_{1,2k+j}]\\ =\min\{ d((-1)^k\e a_{1,2k+1}),\alpha_{2k+j}\}$.
(If $j=m'$, so $2k+j=m$, then $\alpha'_j$ and $\alpha_{2k+j}$ should
be ignored.) But $\alpha'_j=\alpha_{2k+j}$ and in $\fff/\ffs$ we have
$a_{1,2k+j}=(-1)^ka'_{1,j}$ and so
$d((-1)^k\e a_{1,2k+j})=d(\e a'_{1,j})$. It follows that 
$d[(-1)^k\e a_{1,2k+j}]=d[\e a'_{1,j}]$. 

If $1<i\leq j\leq m'$ then
$d[\e a'_{i,j}]=\min\{ d(\e a'_{i,j}),\alpha'_{i-1},\alpha'_j\}$ and
$d[(-1)^k\e a_{2k+i,2k+j}]
=\min\{ d((-1)^k\e a_{2k+i,2k+j}),\alpha_{2k+i-1},\alpha_{2k+j}\}$.
(Again, if $j=m'$, then $\alpha'_j$ and $\alpha_{2k+j}$ should be
ignored.) In $\fff/\ffs$ we have $a'_{i,j}=a'_{1,i-1}a'_{1,j}$,
$a_{2k+i,2k+j}=a_{1,2k+i-1}a_{1,2k+j}$, $a_{1,2k+i-1}=(-1)^ka'_{1,i-1}$
and $a_{1,2k+j}=(-1)^ka'_{1,j}$. It follows that
$a'_{i,j}=a_{2k+i,2k+j}$ and so
$d(\e a'_{i,j})=d(\e a_{2k+i,2k+j})$. Together with
$\alpha'_{i-1}=\alpha_{2k+i-1}$ and $\alpha'_j=\alpha_{2k+j}$, this
implies that $d[\e a'_{i,j}]=d[\e a_{2k+i,2k+j}]$.

(iii) We refer to the proof of Lemma 4.8. If $a_{1,2k}=(-1)^k$ in
$\fff/\ffs$, then $M=J\perp M''$, where
$J=\prec a_1,\ldots,a_{2k}\succ\cong\frac 12A(0,0)^k$ and
$M''=\prec a_{2k+1},\ldots,a_m\succ$. From the cancellation law
[OM, 93:14], since
$M\cong\frac 12A(0,0)^k\perp M'\cong\frac 12A(0,0)^k\perp M''$, we
have $M'\cong M''\cong\prec a_{2k+1},\ldots,a_m\succ$. So the BONG of
$M'$ can be taken such that the sequence $a'_1,\ldots,a'_{m'}$ is
$a_{2k+1},\ldots,a_m$. For every $1\leq i\leq m'$ we have
$[a_1,\ldots,a_{2k}]=FJ\cong H^k$ and
$[a_{2k+1},\ldots,a_{2k+i}]=[a'_1,\ldots,a'_i]$ so
$[a_1,\ldots,a_{2k+i}]=H^k\perp [a'_1,\ldots,a'_i]$.

If $a_{1,2k}=(-1)^k\Delta$ in $\fff/\ffs$, then we proved that
$\prec a_{2k},a_{2k+1}\succ
\cong\prec\Delta a_{2k},\Delta a_{2k+1}\succ$, which implies that also
$[a_{2k},a_{2k+1}]\cong [\Delta a_{2k},\Delta a_{2k+1}]$, and we have
a splitting $M=J\perp M''$, where
$J=\prec a_1,\ldots,\Delta a_{2k}\succ\cong\frac 12A(0,0)^k$ and 
$M''\cong\prec\Delta a_{2k+1},\ldots,a_m\succ$. Again, by the cancelation
law, we have $M'\cong M''$ so the good BONG of $M'$ can be chosen such
that the sequence $a'_1,\ldots,a'_{m'}$ is
$\Delta a_{2k+1},\ldots,a_m$. Since
$[a_{2k},a_{2k+1}]\cong [\Delta a_{2k},\Delta a_{2k+1}]$, for every
$1\leq i\leq m'$ we have $[a_1,\ldots,a_{2k+i}]\cong
[a_1,\ldots,\Delta a_{2k},\Delta a_{2k+1},\ldots,a_{2k+i}]$. But
$[a_1,\ldots,\Delta a_{2k}]\cong FJ\cong H^k$ and
$[\Delta a_{2k+1},\ldots,a_{2k+i}]=[a'_1,\ldots,a'_i]$. Hence, again,
$[a_1,\ldots,a_{2k+i}]=H^k\perp [a'_1,\ldots,a'_i]$. \qed

As an application, we give an explicit necessary condition for
$n$-universality in the case when $n\geq 3$ is odd, which follows from
the Lemma 4.4 in the case when $n-2k=1$. 


\bco If $n\geq 3$ is odd, $M$ is integral and
$M\cong\prec a_1,\ldots,a_m\succ$ relative to some good BONG,
$R_i=R_i(M)$ and $\alpha_i=\alpha_i(M)$, then $M$ represents all
integral lattices over $n$-ary quadratic spaces of Witt index
$\frac{n-1}2$  if and only if $m\geq n+1$, $R_1=R_3=\cdots =R_n=0$,
$R_2=R_4=\cdots =R_{n-1}=-2e$ we have one of the
cases bellow:

\noindent I (a) $\alpha_n=0$ or, equivalently, $R_{n+1}=-2e$

(b) If $m=n+1$ or $R_{n+2}>1$, then
$[a_1,\ldots,a_{n+1}]\cong H^{\frac{n+1}2}$.

(c) If $m\geq n+2$, $R_{n+2}=1$ and either $m=n+2$ or $R_{n+3}>2e+1$,
then $[a_1,\ldots,a_{n+1}]\cong H^{\frac{n+1}2}$.

\noindent II (a) $m\geq n+2$ and $\alpha_n=1$.

(b) If $R_{n+1}=1$ or $R_{n+2}>1$, then $m\geq n+3$ and
$\alpha_{n+2}\leq 2(e-[\frac{R_{n+2}-R_{n+1}}2])-1$.

(c) If $R_{n+1}\leq 0$, $R_{n+2}\leq 1$ and either $m=n+2$ or
$R_{n+3}-R_{n+2}>2e$, then
$H^{\frac{n+1}2}\rep [a_1,\ldots,a_{n+2}]$.
\eco
\pf We write $n=2k+1$ so that $k=\frac{n-1}2$. An $n$-ary quadratic
form $W$ has the Witt index $k$ or $k-1$, corresponding to the case
when $W\cong H^k\perp [b]$ for some $b\in\fff$ or
$W\cong H^{k-1}\perp W'$ with $W'$ ternary anisotropic. Then Lemma 4.4
states that $M$ represents all integral lattices over $n$-ary
quadratic spaces of Witt index $k=\frac{n-1}2$ iff
$M\cong\frac 12A(0,0)^k\perp M'$, where $M'$ is $n-2k=1$-maximal.

A first condition for this to happen is that $M$ splits
$\frac 12A(0,0)^k$, i.e. that it satisfies the conditions (i) and (ii)
of Lemma .... Since $2k=n-1$, condition (i) of Lemma 4.8 states that
$R_1=R_3=\cdots =R_{n-2}=0$ and $R_2=R_4=\cdots =R_{n-1}=-2e$.

Next, assuming that $M$ splits $\frac 12A(0,0)^k$, we write
$M\cong\frac 12A(0,0)^k\perp M'$, where, by the cancellation law
[OM, 93.14], $M'$ is unique up to isometries. Then we are left with
the condition that $M'$ is $1$-universal, which is decided by Theorem
2.1. We write $M'$ as in Lemma 4.9,
$M'\cong\prec a'_1,\ldots,a'_{m'}\succ$, with $R_i(M')=R'_i$ and
$\alpha_i(M')=\alpha'_i$. Since $2k=n-1$, we have $m'=m-n+1$,
$R'_i=R_{n+i-1}$ and $\alpha'_i=\alpha_{n+i-1}$ and, by Lemma
4.9(iii), we may assume that $[a_1,\ldots,a_{n+i-1}]\cong
H^k\perp [a'_1,\ldots,a'_i] =H^{\frac{n-1}2}\perp [a'_1,\ldots,a'_i]$.

Since $m'=n-n+1$ and $R'_1=R_n$, the conditions condition that
$m'\geq 2$ and $R'_1=0$ from Theorem 2.1 applied to $M'$ translate to
$m\geq n+1$ and $R_n=0$. Simlarly, the relations involving the
invariants $R'_i$ and $\alpha'_i$ from the conditions I and II of
Theorem 2.1 translate to the similar relations involving $R_i$ and
$\alpha_i$ from the condition I and II of our corollary, via the
relations $R'_i=R_{n+i-1}$ and $\alpha'_i=\alpha_{n+i-1}$. And since
$[a_1,\ldots,a_{n+1}]\cong H^{\frac{n-1}2}\perp [a'_1,a'_2]$ and
$[a_1,\ldots,a_{n+2}]\cong H^{\frac{n-1}2}\perp [a'_1,a'_2,a'_3]$, the
conditions $[a'_1,a'_2]\cong H$ and $H\rep [a'_1,a'_2,a'_3]$ from
Theorem 2.1 translate to $[a_1,\ldots, a_{n+1}]\cong H^{\frac{n+1}2}$
and $H^{\frac{n+1}2}\rep [a_1,\ldots,a_{n+2}]$ respectively.

Finally, since we alsready have the necessity of $R_n=0$, i.e. of
$R_{2k+1}=0$, the condition (ii) of Lemma 4.8 is superfluous. \qed 

\section*{References}

[B1] C.N. Beli, Integral spinor norms over dyadic local fields,
J. Number Theory 102 (2003) 125-182.
\medskip

\noindent [B2] C.N. Beli, Representations of integral quadratic forms
over dyadic local fields, Electronic Research Announcements of the
American Mathematical Society 12, 100-112, electronic only (2006).
\medskip

\noindent [B3] C.N. Beli, A new approach to classification of integral
quadratic forms over dyadic local fields, Transactions of the American
Mathematical Society 362 (2010), 1599-1617.
\medskip

\noindent [B4] C. N. Beli, Representations of quadratic lattices over
dyadic local fields, (2019), preprint. (arXiv:1905.04552)
\medskip

\noindent [OM] O. T. O’Meara, Introduction to Quadratic Forms,
Springer-Verlag, Berlin (1963).
\medskip

\noindent [XZ] Xu Fei and Zhang Yang, On indefinite and potentially
universal quadratic forms over number fields,
preprint. (arXiv:2004.02090)

\end{document}